\documentclass[10pt]{amsart}

\usepackage{latexsym}
\usepackage{amssymb}
\usepackage{graphicx}
\usepackage{latexsym}
\usepackage{amssymb}

\newcommand{\excise}[1]{}

\numberwithin{equation}{section}

\newtheorem{thm}{Theorem}[section]
\newtheorem{lemma}[thm]{Lemma}

\newtheorem{cor}[thm]{Corollary}

\newtheorem{theorem}{Theorem}

\newtheorem{Example}[thm]{Example}
\newtheorem{Remark}[thm]{Remark}
\newtheorem{Alg}[thm]{Algorithm}
\newtheorem{Defn}[thm]{Definition}
\newtheorem{formula}[thm]{Formula}

\newenvironment{remark}{\begin{Remark}\rm}
                {\mbox{}~\hfill$\square$\end{Remark}}

\newenvironment{proofof}[1]{\begin{trivlist}\item {\it
        Proof of {#1}.\,}}{\mbox{}\hfill$\square$\end{trivlist}}

    {\begin{eqnarray}}
    {\end{eqnarray}$\!\!$}

\newenvironment{eq*}%
    {\begin{eqnarray*}}
    {\end{eqnarray*}$\!\!$}

    {\begin{equation}}
    {\end{equation}$\!\!$}

\newenvironment{eqn*}%
    {\begin{equation*}}
    {\end{equation*}$\!\!$}

        {\begin{list}
                {\noindent\makebox[0mm][r]{\arabic{enumi}.}}
                {\leftmargin=5.5ex \usecounter{enumi}}
        }
        {\end{list}}

\newenvironment{romanlist}%
        {\begin{list}
                {\noindent\makebox[0mm][r]{(\roman{enumi})}}
                {\leftmargin=5.5ex \usecounter{enumi}}
        }
        {\end{list}}

\def\<{\langle}
\def\>{\rangle}
\def\0{{\mathbf 0}}
\def\1{{\mathbf 1}}

\def\CC{{\mathbb C}}

\def\HH{{\mathcal H}{}}

\def\PP{{\mathbb P}}

\def\RR{{\mathbb R}}

\def\TT{{\mathcal T}}
\def\ZZ{{\mathbb Z}}

\def\e{{\varepsilon}}
\def\II{{\mathcal I}}
\def\fib{{\mathbf P}}

\def\th{{\rm th}}

\def\sgn{{\rm sgn}}
\def\pt{{\rm pt}}
\def\tCP{{\widetilde{\CC\PP}}}

\def\pr{{\rm \rho}}
\def\ind{{\rm ind}}

\begin{document}

\title[Morse theory on Hamiltonian $G$-spaces and equivariant K-theory]%
    {Morse theory on Hamiltonian $G$-spaces and equivariant K-theory}

\author{Victor Guillemin}

\author{Mikhail Kogan}

\address{Massachusetts Institute of Technology\\Cambridge, MA}
\email{vwg@math.mit.edu}

\address{Institute for Advanced Study\\Princeton, NJ}
\email{mish@ias.edu}
\thanks{MK was supported by the NSF Postdoctoral Fellowship}


\begin{abstract} Let $G$ be a torus and $M$ a compact Hamiltonian
$G$-manifold with finite fixed point set~$M^G$. If $T$ is a circle
subgroup of $G$ with $M^G=M^T$, the $T$-moment map is a Morse
function. We will show that the associated Morse stratification of
$M$ by unstable manifolds gives one a canonical basis of $K_G(M)$.
A key ingredient in our proof is the notion of local index
$I_p(a)$ for  $a\in K_G(M)$ and $p\in M^G$. We will show that
corresponding to this stratification there is a basis~$\tau_p$,
$p\in M^G$, for $K_G(M)$ as a module over $K_G(\pt)$ characterized
by the property: $I_q\tau_p=\delta^q_p$. For $M$ a GKM manifold we
give an explicit construction of these $\tau_p$'s in terms of the
associated GKM graph.
\end{abstract}

\maketitle

\tableofcontents

\section{Introduction.}
\label{sec:intro} Let $M^{2d}$ be a compact symplectic manifold,
$G$ an $n$-dimensional torus and $\sigma:G\times M\to M$ a
Hamiltonian action. Assume the fixed point set $M^G$ is
finite. Let $T$ be a circle subgroup of~$G$ with the property that
$M^T=M^G$ and let $\phi:M\to \RR$ be the $T$ moment map. This
function is a Morse function and all its critical points are of
even index; so, by standard Morse theory, the unstable manifolds
of $\phi$ with respect to a $G$-invariant Riemannian metric define
a basis of $H_*(M,\RR)$ and by Poincare duality a basis for
$H^*(M,\RR)$ consisting of the Thom classes of the closures of
unstable manifolds. Moreover, these unstable manifolds are
$G$-invariant so they also define a basis for $H_G^*(M)$ as a
module over $H^*_G(\pt)$.

In K-theory the situation is a little more complicated. The critical
points of $\phi$
carry a natural partial order, which is defined by setting $p\leq q$ if
$q$ is inside the closure of the unstable manifold of $\phi$ at $p$ and
then completing this order by transitivity. So, for any unstable
manifold $U$ of $\phi$ at $p$ one can consider the union
$$
W_U=\bigcup U_q
$$
of unstable manifolds $U_q$ for $q\geq p$. It is known that the
there exist classes in $K$-theory which are supported on this set.
However, except in certain special cases (e.g. algebraic torus
actions), it is not known whether there is a genuine (Thom) class
in K-theory associated with $U$. (For algebraic torus actions such
classes can be defined using the structure sheaf of the closure of
$U$, see~\cite{BFM} for details).

We will show in this paper, however, that there is another way of
attaching to the Morse decomposition of~$M$ a basis of~$K_G(M)$
which works even in the case of nonalgebraic torus actions. (As
will be explained later, in the algebraic case our classes will be
different from those constructed using structure sheaves.) The key
idea in our approach is a notion of \emph{local index} for a
K-class $a\in K_G(M)$ at a critical point p of $\phi$. This is
defined as follows: Let $S$ be the stable manifold of~$\phi$
at~$p$, and for small $\e>0$ let $S_\varepsilon$ be the compact
symplectic orbifold obtained from $S$ by the symplectic cutting
operation of Lerman~\cite{Lerman}. We recall that $S_\varepsilon$
is obtained from the manifold with boundary
\begin{equation}
\label{eq:cut} \tilde S_\varepsilon =\{x\in S,\phi(x)\geq
\phi(p)-\varepsilon\}
\end{equation}
by collapsing to points the $T$-orbits on the boundary. In
particular, there is a projection $\pr: \tilde S_\e\to S_\e$ and
an inclusion $i:\tilde S_\e\to M$; so a K-class $a\in K_G(M)$
defines a class $\kappa_\e(a)=\pr_!i^*a\in K_G(S_\e)$, where
$\pr_!$ is the pushforward map (we will define the map $\kappa_\e$
in more detail in Section~\ref{sec:local-index}).

Now let \emph{the local index} of $a$ at $p$
$$
I_p(a)\in K_G(\pt)
$$
be the Atiyah-Segal index of $\kappa_\e(a)$, that is, the
pushforward of $\kappa_\e(a)$ with respect to the map $S_\e\to
\pt$. Recall that $K_G(\pt)$ is just the representation ring
$R(G)$ of the torus $G$, so that each local index is just a
virtual representation of $G$. One of the main results of this
paper is the following theorem.

\begin{theorem}
\label{thm:main} Let $p$ be a critical point of $\phi$ and $U$ the
unstable manifold of $\phi$ at $p$. Then there exists a unique
K-theory class $\tau_p\in K_G(M)$ with the properties:
\begin{romanlist}
\item
\label{cond1} $I_p(\tau_p)=1$,
\item
\label{cond2} $I_q(\tau_p)=0$ for all critical points $q$ of
$\phi$ except $p$,
\item
\label{cond3} The
restriction of $\tau_p$ to a critical point $q$ is zero unless
$q\in W_U$.
\end{romanlist}
Moreover, the  $\tau_p$'s generate $K_G(M)$ freely as a
module over $K_G(\pt)$.
\end{theorem}

Let $\II:K_G(M)\to K_G(M^G)$ be the map which takes the value
$I_p$, at $p$.  This we will call the \emph{total index map}. As explained 
in Remark~\ref{rem:explanation} the total index is not an $R(G)$-module 
homomorphism but it is a homomorphism with respect to the subring,
$R(G/T)$, of $R(G)$.  
Theorem~\ref{thm:main} implies

\begin{cor}
\label{cor:main}The total index map, $\II$, is an $R(G/T)$ module
isomorphism.
\end{cor}

\begin{remark}
Notice that we can define local indices even if the fixed point
set of the action is not finite. Namely, let $F$ be a connected
component of $M^G$, not necessarily consisting of one point.
Then if $S$ is the stable manifold of $\phi$ at $F$ we can still
define $S_\e$, the projection $\pr$, the inclusion $i$, and the map
$\kappa_\e$. Moreover, there is a fibration $\fib:S_\e\to F$ whose
fibers are weighted projective spaces. So, we can define the local
index at $F$
$$
I_F: K_G(M)\to K_G(F)
$$
to be $\fib_!\kappa_\e$, the composition of the pushforward $\fib_!$ with 
$\kappa_\e$. Then the total local index map
$$
\II: K_G(M)\to K_G(M^G)
$$
is well-defined and is an $R(G/T)$-module homomorphism; so, it is
natural to pose the following question whose answer, we believe,
depends on whether or not $\phi$ is K-theoretically perfect.

\medskip

\noindent{\bf Question.} \it When is $\II$ an isomorphism?
\end{remark}

\begin{remark} Notice that local indices can also be defined in the 
setting of equivariant cohomology. Namely, for $a\in H_G^*(M)$, we 
let $I_p(a)$ be the pushforward (or integral)
of~$\kappa'_\e(a)$, where $\kappa'_\e:H^*_G(M)\to H^*_G(S_\e)$ is defined 
in exactly the same way as $\kappa_\e$. It will be clear from the proof of 
Theorem~\ref{thm:main} that its analogue for equivariant
cohomology is also true, and that the cohomological analogues of
the $\tau_p$'s  are just ``the equivariant 
Poincare duals'' of the closures of the unstable manifolds.
\end{remark}

The other main result of this paper is a constructive version of
Theorem~\ref{thm:main} for GKM spaces, that is, an explicit
computation of the classes $\tau_p$. We start by recalling some
facts about GKM spaces. The \emph{one-skeleton} of $M$
\begin{equation}
\label{eq:fixed} \{x\in M, \ \dim G\cdot x=1\}
\end{equation}
is a union of symplectic submanifolds of $M$. The action $\sigma$
is defined to be \emph{a GKM action} and $M$ \emph{a~GKM space} if
each connected component of the one-skeleton is exactly of
dimension $2$. It is easy to see that if $\sigma $ is GKM the
fixed point set $M^G$ has to be finite. Let
$$
V=\{p_1,\dots, p_\ell\}
$$
be the points of $M^G$. To each connected component $e_i^\circ$
of~(\ref{eq:fixed}) let $e_i$ be its closure, and let
$$
E=\{e_1,\dots, e_N\}
$$
be the set of $e_i$'s. We claim:
\begin{romanlist}
\item $e_i$ is an imbedded copy of $\CC\PP^1$,
\item $e_i-e_i^\circ$ is a two element subset of $V$,
\item for $i\neq j$ the intersection $e_i\cap e_j$ is empty or is a
one-element subset of $V$,
\item for $p\in V$ the set $\{e_i\in E, p\in e_i\}$ is $d$-element
subset of $V$.
\end{romanlist}
For the proof of these assertions see for instance \cite{GZ1}.
These assertions can be interpreted as saying that $V$ and $E$ are
the vertices and edges of a $d$-valent graph $\Gamma$.

One can describe the action of $G$ on the
one-skeleton~(\ref{eq:fixed}) by means of a labeling function
which labels each oriented edge of this graph by an element of the
weight lattice $\ZZ^*_G$ of $G$. Explicitly, let $e$ be an edge of
$\Gamma$ joining a vertex $p$ and a vertex $q$. To $e$ we can
associate two oriented edges $e_p$ and $e_q$ pointing from $p$ to
$q$ and from $q$ to $p$ respectively.  As a geometric object $e$
is a $G$ invariant imbedded $\CC\PP^1$ with fixed points at $p$
and $q$; if we denote by $\alpha_{e_p}$ the weight of the isotropy
representation of $G$ on the tangent space to $e$ at $p$ (and by
$\alpha_{e_q}=-\alpha_{e_p}$ the weight of the isotropy
representation at~$q$) we get a labeling function $\alpha$ which
describes how each connected component of the one-skeleton is
rotated about its axis of symmetry by~$G$.

GKM theory is concerned with reconstructing, in so far as
possible, the geometry of $M$ from the combinatorics of the pair
$(\Gamma,\alpha)$. It is known for instance that the ring
structure of~$H_G^*(M)$ and~$K_G(M)$ are determined by $(\Gamma,
\alpha)$. (See \cite{GKM, At, CS, KR, TW} for versions of this
result. We'll explain below how $K_G(M)$ is determined by
$(\Gamma, \alpha).$) It was also shown in~\cite{GZ3} that if $T$
is generic circle subgroup of $G$, and $\tau\in H_G^*(M)$ the
cohomology class dual to an unstable manifold of $T$, the
restriction of $\tau$ to $M^G$ is completely determined by
$(\Gamma, \alpha)$. In this paper we will prove analogous results
for $K_G(M)$.

Let us recall how the ring structure of $K_G(M)$ is determined by
$(\Gamma, \alpha)$. One knows that the restriction map
\begin{equation}
\label{eq:restriction} K_G(M)\to K_G(M^G)
\end{equation}
is an injection, so $K_G(M)$ is a subring of the much simpler ring
\begin{equation}
\label{eq:subring} K_G(M^G)=\bigoplus_{i=1}^\ell K_G(p_i)
\end{equation}
Since $K_G(\pt)=R(G)$, an element of the ring~(\ref{eq:subring})
is just a map
\begin{equation}
\label{eq:map} \chi:V\to R(G)
\end{equation}
and one has:

\begin{theorem} {\rm{\cite{At,KR}}}
For each $e\in E$ let $G_e$ be the kernel of the homomorphism
$$
e^{2\pi \sqrt{-1}\alpha_e}:G\to S^1.
$$
Then the element~(\ref{eq:map}) of $K_G(M^G)$ is in the image
of~(\ref{eq:restriction}) if and only if for every $e\in E$.
\begin{equation}
\label{eq:edge} r_e(\chi_p)=r_e(\chi_q)
\end{equation}
$p$ and $q$ being the vertices of $e$ and $r_e$ the restriction
map $R(G)\to R(G_e)$.
\end{theorem}

For GKM manifolds one can also translate some aspects of
Morse theory into the language of graphs. Recall that $\phi$ is
the moment map on $M$ with respect to the circle $T$ action. Think
of each edge $e$ of the graph connecting vertices $p$ and $q$ as
two oriented edges $e_p$ and $e_q$. Then if $\phi(p)>\phi(q)$ we
say that the edge $e_p$ going from $p$ to $q$ is \emph{descending}
and $e_q$ from $q$ to $p$ \emph{ascending}. If $U$ is the unstable
manifold of $\phi$ at $p$ then every fixed point, $q$, inside $W_U$
is the terminal point of a path on $\Gamma$ starting at $p$ and
consisting of ascending edges; and this gives one a way of
describing $W_U$ in terms of $\Gamma$.  In particular, we will  present below an explicit formula for the image of
$\tau_p$ under the imbedding~(\ref{eq:restriction}), which
expresses the restriction of $\tau_p$ to $q\in M^G$ as a sum of
combinatorial expressions associated with the ascending paths in
$\Gamma$ going from $p$ to $q$. (An analogous formula for the
cohomological counterpart of $\tau_p$ can be found in~\cite{GZ3}.)
Our results will follow from the following theorem, which allows one
to compute local indices  in terms of restrictions of K-theory
classes to fixed points and vice versa.

\begin{theorem}
\label{thm:GKM} For $p\in V=M^G$, let  $e_1,\dots, e_m$ be
the descending edges with initial vertex at $p$. Let the edge $e_i$ connect $p$
to $q_i$ and be labeled by the weight $\alpha_i$. Then for any $a\in
K_G(M)$ we have
\begin{equation}
\label{eq:GKM} I_p(a)=\sum_{i=1}^m \tilde \pi_i \tilde
r_i\Big(\frac{a_{q_i}}{(1-\zeta)\prod_{j\neq i} \big(1-e^{2\pi
\sqrt{-1} \alpha_{j}}\big)}\Big) + \frac{a_p}{\prod_{i=1}^m
(1-e^{2\pi \sqrt{-1}\alpha_{i}})},
\end{equation} 
where $a_q$ is the restriction of $a$ to $q$, $\zeta$ is the generator of
the character
ring $R(T)$, $\tilde r_i$ is  the
restriction $R(G\times T)\to R(G_{e_i}\times T)$ and $\tilde
\pi_i:R(G_{e_i}\times T)\to R(G)$ is the Gysin map defined
in~(\ref{eq:pi2}).
\end{theorem}

We conclude this introduction with a section-by-section summary of
the contents of the paper. In Section~\ref{sec:proofs} we prove
Theorem~\ref{thm:main} by adopting certain arguments from 
classical Morse theory to the setting of equivariant K-theory,
and in Section~\ref{sec:lagrance} we prove an analogue for the ring
$R(G)$ of the ``Lagrange interpolation formula'' 
 of~\cite{GZ3}. More specifically,
let $T$ be a circle subtorus of $G$ and $H\subset G$  a
complimentary subtorus, so that $G=T\times H$. Let $\hat R(G)$ be
the ring of polynomials $\sum_kc_kz^k$ where $c_k$ is in the
quotient ring $Q(H)$ of $R(H)$. Let $w:T\to S^1$ be an
isomorphism. Via the splitting $G=T\times H$ we can extend $w$ to
a homomorphism of~$G$ onto $S^1$ by setting~$w$ equal to $1$ on
$H$. Let $\xi$ be the infinitesimal generator of $T$ (chosen so
that it corresponds under $w$ to the standard generator
$\frac{\partial}{\partial \theta}$ of $S^1$). The interpolation in
question is with respect to weights $\alpha_i\in \ZZ^*_G$,
$i=1,\dots, m$.  Letting~$G_i$ be the kernel of the homomorphism
$$
e^{2\pi\sqrt{-1}\alpha_i}: G\to S^1
$$
it describes to what extend an element $f$ of $R(G)$ is determined
by its restrictions to the $R(G_i)$'s. More explicitly it asserts:

\begin{theorem}
\label{thm:lagrange} If the $\alpha_i$'s are pairwise linearly
independent and $\alpha_i(\xi)\neq 0$, then there exist elements
$f_0\in R(G)$ and $f_i\in\hat R(G)$, $i=1,\dots, m$ such that
\begin{equation}
\label{eq:lagrande-intro}
\frac{f}{\prod_{i=1}^m(1-e^{2\pi\sqrt{-1}
\alpha_i})}=f_0+\sum_{i=1}^m \sgn(\alpha_i(\xi)) f_i
\end{equation}
where
\begin{equation}
\label{eq:lagrange-intro2} f_i(z)=\pi_i r_i
\Big(\frac{f(w,h)}{(1-\frac{z}{w})\prod_{j\neq
i}(1-e^{2\pi\sqrt{-1} \alpha_j})}\Big).
\end{equation}
Here $r_i$ is the restriction map $R(G)\to
R(G_i)$, and $\pi_i$ is the Gysin map $R(G_i)\to R(H)$ associated
with the projection $G_i\to G/T\cong H$ (see (\ref{eq:Gysin}) for
definitions).
\end{theorem}

In Section~\ref{sec:index} we will apply the Atiyah-Segal
localization theorem to equivariant $K$-classes on twisted
projective space to obtain a formula for the equivariant index of
such a class, and in  Section~\ref{sec:local-index} we will apply this formula to the
 twisted projective space, $S_\e$, and show that the
formulas~(\ref{eq:GKM}) and (\ref{eq:lagrande-intro}) are essentially the 
same formula viewed from different perspectives, i.e.,~are the 
topological and algebraic versions of this formula. More specifically, we 
will show that if $p$ is a vertex of the GKM graph $\Gamma$, $f$ the 
restriction to~$p$
of an element $a$ of $K_G(M)$ and the $\alpha_i$'s the weights
associated with the descending edges of~$\Gamma$ with initial vertex at $p$
then the $f_0$ in~(\ref{eq:lagrande-intro}) is just the local
index $I_p(a)$. We will then use this result to prove
Theorem~\ref{thm:GKM}.

In Section~\ref{sec:restriction} we will obtain explicit
formulas for the $\tau_p$'s in terms of their restrictions to the
fixed points.   These will be proved by a repeated iteration of
(1.7).  (We recall that $I_q (\tau_p) =0$ if $p \neq q$ and $I_p
(\tau_p)=1$; so (1.7) gives one an effective way of computing
$\tau_p$ at $q$ in terms of the values of $\tau_p$ at the points
in $M^G$ lying below $q$ in $W_U$.)  In the last section we give a more general
definition of the notion of ``local index'', for which many of the results
 above still hold with  minor modifications.  This new definition
 involves choosing, for each $p \in M^G$, a circle subgroup,
 $T_p$, depending on $p$, and replacing the space,
 $S_{\epsilon}$, by the space obtained by cutting the stable
 manifold, $S$ of $\phi$ at $p$ by $T_p$. If these spaces are manifolds,
i.e.,~don't
have orbifold singularities, the formulas (\ref{eq:GKM}) and 
(\ref{eq:lagrande-intro}) become
considerably simpler (for instance the Gysin maps in these
formulas are all identity maps).  In particular these formulas
are now very similar to the analogous formulas in
equivariant cohomology.  (See \cite{GZ3}.)  As an application of these results, we discuss this generalized
index map for the Grassmannian and explain some tie-ins of our
results with recent work of Lenart \cite{Lenart} on Shur and Grothendieck 
polynomials.

\section{Morse theory and equivariant K-theory.}
\label{sec:proofs}

We will deduce Theorem~1 from the
series of lemmas below.  These lemmas are K-theoretic analogs of classical results in equivariant
Morse theory~\cite{AB}.

As before, $M$ is a compact symplectic manifold with a Hamiltonian
$G$ action, $T$ is a generic circle subgroup of~$G$ with $M^T=M^G$
and $\phi$ is the moment map of the $T$ action. For a critical point
$p$ of $\phi$, that is, $p\in M^G$, pick a $G$-invariant complex
structure on the tangent space $T_p$ at $p$ compatible with the
symplectic structure. Then  $T_p$  splits into the negative and
positive components $T_p^-$ and~$T_p^+$ on which the circle $T$
acts with negative and positive weights respectively. Let
$\Lambda^-_p$ be the 
virtual vector space $\sum (-1)^k \Lambda^k(T_p^-)$ with its given
$G$ action. By definition this
is an element of $R(G)\cong K_G(p)$. Moreover, if $\alpha_1,\dots,
\alpha_m$ are the (possibly repeating) weights of the $G$ action on
$T^-_p$, then,  as a virtual character of
$G$
\begin{equation}
\label{eq:lambda} \Lambda_p^-
=\prod_i(1-e^{2\pi\sqrt{-1}\alpha_i}).
\end{equation}

Recall that $K_G^n(M)$ is the compactly supported
K-group, $K_{G,c}(M\times \RR^n)$ and that $K_G(M)=K^0_G(M)$.
Moreover, by 
Bott periodicity
$$
K_G^n(M)\cong K_G^{n+2}(M).
$$

For a critical point $p$, let $\phi(p)=c$. Without loss of
generality we may assume there is only one critical point $p$ in
$\phi^{-1}(c)$. For a small $\e>0$, let
$$
M_p^+=\{x\in M|\phi(x)\leq c+\e\}\text{ and }M_p^-=\{x\in
M|\phi(x)\leq c-\e\}.
$$

\begin{lemma} \label{lem:sequence}
The K-theory long exact sequence  for the pair $(M_p^+, M_p^-)$
splits into short exact sequences
\begin{equation}
\label{eq:sequence} 0\to K^*_G(M_p^+,M_p^-)\to K^*_G(M_p^+)\to
K^*_G(M_p^-)\to 0.
\end{equation}
\end{lemma}

\begin{proof}\footnote{We thank Sue Tolman and Jonathan Weitsman for suggesting
this argument to us.} Let $S_p$ be the stable manifold at $p$ and
$S^-_p=S_p\cap M_p^-$. Then there are isomorphisms
$$
K^*_G(p)\stackrel{\TT}{\to} K^*_G(S, S^-) \stackrel{\HH}{\to}
K^*_G(M_p^+,M_p^-)
$$
where $\TT$ is the Thom isomorphism and $\HH$ comes from homotopy
equivalence.

To show that the long exact sequence splits it is enough to show
that the maps
$$
\mathcal J: K^*_G(M_p^+,M_p^-)\to K^*_G(M_p^+)
$$
are all injective.

Let $\iota_p$ be the inclusion of $p$ into $M^+$. It is well known
that for  $*=0$, the map
$$
\iota_p^*\circ \mathcal J\circ \HH\circ \TT:K^0_G(p)\to K^0_G(p)
$$
is a multiplication by $\Lambda^-_p$ and hence injective, since
by~(\ref{eq:lambda}) $\Lambda^-_p$ is not a zero divisor in
$R(G)$. Therefore the map $\mathcal J$ must also be injective for
$*=0$.

For $*=1$, notice that by the Thom isomorphism
$$
K^1_G(S, S^-)=K^0_{G,c}(\RR)
$$
where the action of $G$ on $\RR$ is trivial. Since the
nonequivariant K-space, $K^1(\pt)=K^0_c(\RR)$ is known to be
trivial, it is easy to conclude that $K^0_{G,c}(\RR)$ is trivial
as well. Indeed, by splitting each vector bundle on $\RR$ into
bundles for which the action  of $G$ on the fibers is given by a
single weight, we reduce the calculation of $K^0_{G,c}(\RR)$ to
the calculation of $K^0_c(\RR)$. This implies that $K^1_G(S, S^-)$
is trivial, and that $\mathcal J$ is injective, and hence finishes the
proof.
\end{proof}

\begin{cor}
\label{cor:restriction} The restriction
$$
K_G(M)\to K_G(M^G)
$$
is injective.
\end{cor}

\begin{proof}
For every critical point $p$ it is enough to show that if the
restriction
\begin{equation}
\label{eq:inject1} K_G(M_p^-)\to K_G(M_p^-\cap M^G)
\end{equation}
is injective then
\begin{equation}
\label{eq:inject2} K_G(M_p^+)\to K_G(M_p^+\cap M^G)
\end{equation}
is also injective.

The long exact sequence of the pair $(M_p^+\cap M^G,M_p^-\cap
M^G)$ obviously splits into short exact sequences. By
Lemma~\ref{lem:sequence} the short exact
sequence~(\ref{eq:sequence}) maps into the corresponding short
exact sequence for $(M_p^+\cap M^G,M_p^-\cap M^G)$. The
restriction
$$
K_G(M_G^+, M_G^-)\to K_G(p)
$$
is an injection, since $\Lambda^-_p$ is not a zero divisor. So,
the injectivity of~(\ref{eq:inject1}) together with the five-lemma
implies the injectivity of~(\ref{eq:inject2}).
\end{proof}

Recall that the elements of $M^G$ are partially ordered with $q\leq p$ if $q$ lies in the
closure of the unstable manifold at $p$. Completing this relation by
transitivity one gets a partial order on $M^G$. The proof of
the following Lemma is postponed until
Section~\ref{sec:local-index} where we will deduce it from  the Atiyah-Segal
localization theorem for weighted projective spaces.

\begin{lemma}
\label{lem:local} Assume $\tau\in K_G(M)$ restricts to zero at
every $q\in M^G$ with $q<p$. Then
$$
\tau(p)=I_p(\tau) \Lambda^-_p
$$
where $\tau(p)$ is the restriction of $\tau$ to $p$.
\end{lemma}

\begin{lemma}
\label{lem:existence} For every critical point $p$ of $\phi$ there
exists an element $\tau$ of $K_G(M)$ which restricts to zero at
every $q\in M^G$  with $\phi(q)<\phi(p)$ (in other words, $\tau$
is supported above~$p$), such that $I_p(\tau)=1$.
\end{lemma}

\begin{proof}
By Lemma~\ref{lem:local}, it is enough to construct a K-class, $\tau_p$,
supported above $p$ whose restriction to the point~$p$
is~$\Lambda^-_p$. By Lemma~\ref{lem:sequence} such an element
exists in $K_G(M^+_p)$. Moreover, by induction such an element
 exists in $K_G(M^+_q)$ for every $q$ with $\phi(q)>\phi(p)$.
Indeed, because of the short exact sequence~(\ref{eq:sequence})
for the pair $(M^+_q,M^-_q)$ if such an element exists in
$K_G(M_q^-)$ it can be lifted to an element in $K_G(M_q^+)$.
\end{proof}

\begin{proofof}{Theorem~\ref{thm:main}}

Let $p_1,\dots, p_r$ be  the  points of $M^G$ ordered such that
$\phi(p_1)\geq\phi(p_2)\geq\dots\geq \phi(p_r)$. Let us prove by induction
on $k$ that classes $\tau_{p_k}$ satisfying (\ref{cond1}), (\ref{cond2})  
and (\ref{cond3}) exist. For $k=1$ this follows from
Lemmas~\ref{lem:local} and~\ref{lem:existence}.  Assume we can construct
classes $\tau_{p_1},\dots, \tau_{p_{k-1}}$ satisfying these properties. By
Lemma~\ref{lem:existence} we can choose a class $\tau$ supported above $p$
with $I_{p_k}(\tau)=1$. By Lemma~\ref{lem:local} we conclude that
$I_{p_\ell}(\tau)=0$ for every $\ell>k$. Choose the largest $\ell<k$ with
$I_{p_\ell}(\tau)\neq 0$. Then change $\tau$ to
$\tau-I_{p_\ell}(\tau)\tau_{p_\ell}$. It is clear that the new $\tau$
satisfies $I_{p_\ell}(\tau)=\delta^\ell_k$  for
$m\geq \ell$. Hence by induction we can find a class $\tau_{p_k}$ such
that $I_{p_\ell}(\tau_{p_k})=\delta^m_k$ for all $m$, implying
(\ref{cond1})  and (\ref{cond2}).


To show that such $\tau_{p_k}$ satisfies~(\ref{cond3}), assume that for 
every $m<\ell$ either $p_m\geq p_k$, or $\tau_{p_k}(p_m)=0$. 
We will now show that if $p_\ell\ngeq p_k$ then $\tau_{p_k}(p_\ell)=0$.
Notice that every $p<p_\ell$ is among the points 
$p_{\ell+1},\dots, p_r$ satisfying $p\ngeq p_k$ so that $\tau_{p_k}(p)=0$ 
by the above assumption. 
Hence $\tau_{p_k}(p_{\ell})=0$ by
Lemma~\ref{lem:local}.


To prove uniqueness of the classes $\tau_p$ it remains to show that if
both $\tau_p$ and $\tau_p'$ satisfy (\ref{cond1})-(\ref{cond3}),
then they must be equal. Set $\delta=\tau_p-\tau_p'$. Then
$\delta$ is a K-theory class whose local indices are all zero. By
Corollary~\ref{cor:restriction}, if~$\delta$ were not zero, there
would exist a critical point~$q$ such that the restriction
of~$\delta$ to~$q$ is nonzero. Pick such $q$ with
minimal~$\phi(q)$. Then by Lemma~\ref{lem:local} the local index
$I_q(\delta)$ is not zero. 

Finally, we need to show that $\tau_p$'s generate $K_G(M)$ freely as a
$K_G(\pt)$ module.  We first prove that every $a\in K_G(M)$ can be 
expressed as linear a combination of $\tau_p$'s with coefficients 
in $R(G)$. 
To do so order the points 
of $M^G$ as above. We will prove by induction on $k$ that if $a$ 
restricts to zero at $p_{k+1},\dots, p_r$, then $a$ is a linear 
combination of $\tau_{p_1},\dots, \tau_{p_k}$. Clearly this holds for
$k=1$. To prove the induction step for $k>1$, 
notice that if $a$ is supported above $p_k$ the class
$a-I_{p_k}(a)\tau_{p_k}$ restricts to zero at 
$p_{k},\dots, p_r$ and by the induction assumption is a linear combination of 
$\tau_{p_1},\dots, \tau_{p_{k-1}}$. This proves the induction step. It 
remains to show that no nontrivial linear
combination, $\gamma=\sum c_{k}\tau_{p_k}$, with $c_k\in R(G)$ is zero.
Because of Corollary~\ref{cor:restriction} it is enough to show
that $\gamma$ restricts nontrivially to $M^G$ whenever one of
$c_k$'s in not zero. Pick the largest $k$ with $c_k\neq 0$. 
Then by Lemma~\ref{lem:local}
the restriction of $\delta$ to $p_k$ is not zero. This finishes the
proof of Theorem~\ref{thm:main}.
\end{proofof}

\begin{proofof}{Corollary~\ref{cor:main}} The injectivity of $\II$ follows 
from the fact that  $a\in K_G(M)$ must be equal to $a'\in K_G(M)$ 
whenever $I_p(\tau)=I_p(\tau')$ for  all $p\in M^G$, which is true, since
by 
the  argument used in the proof of Theorem~\ref{thm:main} the class $a-a'$ 
must vanish.  

To prove surjectivity we  must show that for any choice of $i_p\in R(G)$ 
there exists a class $a\in K_G(M)$ with $I_p(a)=i_p$. Order the points 
of $M^G$ as in the proof of Theorem~\ref{thm:main}. Take $a=\sum_{k=1}^r 
c_k\tau_{p_k}$ with 
$$
c_k=i_{p_k}-I_{p_k}\Big(\sum_{\ell=k+1}^r c_\ell\tau_{p_\ell}\Big).
$$
Then it is easy to see that $I_p(a)=i_p$ for  every $p$.
\end{proofof}

\section{Lagrange interpolation.}
\label{sec:lagrance}

We recall the statement of the classical
Lagrange interpolation formula for the ring of polynomials in one
variable and then describe how to generalize this formula to the representation ring
$R(G)$.

Let $f(z)=\sum_{k=0}^N c_kz^k$ be a polynomial in $z$ with complex
coefficients. Given $d$~distinct complex numbers, $a_i\in \CC$,
$i=1,\dots,d$ one has the identity

\begin{formula}
\begin{equation}
\label{eq:lagrange1}
\frac{f(z)}{\prod_i(z-a_i)}=f_0(z)+\sum_i\frac{c_i}{z-a_i}
\end{equation}
where $f_0(z)$ is a polynomial of degree $N-d$ and
\begin{equation}
c_i=\frac{f(a_i)}{\prod_{j\neq i}(a_i-a_j)}.
\end{equation}
\end{formula}

\begin{proof}
The function
$$
g(z)=\frac{\sum_i f(a_i)\prod_{j\neq i}(z-a_j)}{\prod_{j\neq
i}(a_i-a_j)}
$$
is a polynomial of degree $d-1$ which takes the same values as
$f(z)$ at the points $z=a_i$, so $f(z)-g(z)$ is divisible by
$\prod_i(x-a_i)$.
\end{proof}

A slightly more complicated variant of this identity is the
following. Let $k_1,\dots, k_d$ be positive integers. Define the
``Gysin map'' $\pi_i$ which acts on polynomials in two variable
$z$ and $w$ by
\begin{equation}
\label{eq:Gysin} \pi_i(h(z,w))=\frac{1}{k_i}\sum_{\ell=1}^{k_i}
h(z,w_{i,\ell})
\end{equation}
summed over the roots $w_{i,\ell}$ of $z^{k_i}-a_i$. Then

\begin{formula}
\begin{equation}
\label{eq:lagrange2}
\frac{f(z)}{\prod_i(z^{k_i}-a_i)}=f_0(z)+\sum_i\frac{f_i(z)}{z^{k_i}-a_i}
\end{equation}
where $f_0$ is a polynomial of degree $N-\sum k_i$ and $f_i$ is
the polynomial
\begin{equation}
\label{eq:lagrange3}
\frac{z^{k_i}-a_i}{a_i}\pi_i\Big(\frac{f(w)w}{(z-w)\prod_{j\neq
i}(w^{k_j}-a_j)}\Big).
\end{equation}
Moreover, if $f=\sum_{i=0}^N c_iz^i$ then $f_0=\sum_{i=1}^{N-k}
d_j z^j$, where $k=\sum k_i$ and
\begin{equation}
\label{eq:new} d_j=\sum_{i=k+j}^Nc_i\Big(\sum_{\ell_1k_1+\dots
+\ell_d k_d=i-k-j}a_1^{\ell_1}\dots a_d^{\ell_d}\Big).
\end{equation}
\end{formula}

\begin{proof}
Factoring
$$
z^{k_i}-a_i=\prod_{i=1}^{k_i} (z-w_{i,\ell})
$$
and applying~(\ref{eq:lagrange1}) we get
$$
\frac{f(z)}{\prod_i(z^{k_i}-a_i)}=f_0+\sum_{i=1}^d h_i
$$
where
$$
h_i(z)=\sum_{\ell=1}^{k_i}\frac{f(w_{i,\ell})}{z-w_{i,\ell}} \cdot
\frac{1}{\big(\prod_{m\neq\ell}(w_{i,\ell}-w_{i,m}) \big) \big(
\prod_{j\neq i, 1\leq m\leq k_j}(w_{i,\ell}-w_{j,m}) \big)}.
$$
But
$$
\prod_{m\neq\ell}(w_{i,\ell}-w_{i,m}) =\lim_{w\to w_{i,\ell}}
\frac{w^{k_i}-a_i}{w-w_{i,\ell}}=k_iw^{k_i-1}|_{w_{i,\ell}}=\frac{k_ia_i}{w_{i,\ell}}
$$
and
$$
\prod_{j\neq i, 1\leq m\leq k_j}(w_{i,\ell}-w_{j,m})= \prod_{j\neq
i} (w_{i,\ell}^{k_j}-a_j).
$$
So
$$
h_i(z)=\frac{1}{k_i}\sum_{\ell=1}^{k_i}h_i(z,w_{i,\ell})=\pi_i(h_i(z,w))
$$
where
$$
h_i(z,w)=\frac{1}{a_i}\frac{f(w)w}{(z-w)\prod_{i\neq
j}(w^{k_j}-a_j)}.
$$
This proves~(\ref{eq:lagrange2}) and~(\ref{eq:lagrange3}).

To prove~(\ref{eq:new}) expand both sides of~(\ref{eq:lagrange2})
in  powers of $z$. In particular, note that
\begin{equation}
\label{eq:expansion}
\frac{f(z)}{\prod_i(z^{k_i}-a_i)}=\sum_{i=0}^N c_i
z^{i-k}\prod_{i=1}^d \Big(1-\frac{a_i}{z^{k_i}}\Big)^{-1}=
\sum_{i=0}^Nc_i\Big(\sum_{\ell_1,\dots ,\ell_d=0}^\infty
a_1^{\ell_1}\dots a_d^{\ell_d}\Big)z^{i-k-\sum\ell_ik_i}.
\end{equation}
It is clear from~(\ref{eq:lagrange2}) that $f_0$ in the
``polynomial part'' of this expression. Indeed the expansion of
$\frac{f_i(z)}{z^{k_i}-a_i}$  only involves negative powers of $z$ and
$f_0$ is a polynomial in $z$. Therefore, $f_0$ is the sum of the
terms in the expression~(\ref{eq:expansion}) for which the exponent
$i-k-\sum\ell_ik_i$ is greater than or equal to zero. Hence its
coefficients are given by~(\ref{eq:new}).
\end{proof}

A third variant of this formula involves the character ring of the group
$\CC^*=\CC-\{0\}$, i.e. finite sums of the form
$f(z)=\sum_{i=-N}^Mc_kz^k$ with $c_k\in\CC$. It asserts that if
$k_1,\dots, k_d$ are non-zero integers, then

\begin{formula}
\label{form:lagrange}
\begin{equation}
\label{eq:lagrande4}
\frac{f(z)}{\prod_i(1-a_iz^{k_i})}=f_0(z)-\sum_i \sgn (k_i)
  \frac{ f_i(z)}{1-a_i z^{k_i}}
\end{equation}
where $f_0, f_1, \ldots ,f_d$ are in the character ring of $\CC^*$,  and
\begin{equation}
\label{eq:lagrange5} f_i(z)= (1-a_i z^{k_i}) \pi_i
\Big(\frac{f(w)w}{(z-w)\prod_{j\neq i}(1-a_jw^{k_j})}\Big)
\end{equation}
\end{formula}

\begin{proof} This is easily deduced from (\ref{eq:lagrange2}) and
(\ref{eq:lagrange3}) by setting
\begin{align*}
1-a_iz^{k_i}&=z^{k_i}(z^{-k_i} -a_i)\text{ for }k_i\text{
negative  and }\\
\text{}&=-a_i(z^{k_i}-\frac{1}{a_i})\text{ for }k_i\text{
positive}
\end{align*}
and applying~(\ref{eq:lagrange2}) to the function
$$
g=f\prod_{k_i>0}\Big(\frac{-1}{a_i}\Big) \prod_{k_i<0} z^{-k_i}
$$
instead of $f$.
\end{proof}

Our last version of Lagrange interpolation is a multidimensional
generalization of the character formula~(\ref{eq:lagrande4}). As
in Theorem~\ref{thm:lagrange} let $G$ be an $n$-dimensional torus,
$R(G)$ the character ring of $G$ and $\alpha_i$, $i=1,\dots,d$
elements of the weight lattice of $G$.
Corresponding to each~$\alpha_i$ one has a homomorphism
$$
e^{2\pi\sqrt{-1} \alpha_i}:G\to S^1.
$$
Let $G_i$ be the kernel of this homomorphism and let $T$ be a
circle subgroup of $G$. Assume $\xi$ is the infinitesimal
generator of $T$. Fixing a complimentary subtorus $H$ to $T$ in
$G$ we have
$$
G=T\times H=S^1\times H\subseteq \CC^*\times H.
$$
Hence one can regard an elements of $R(G)$ as a finite sum $f=\sum
_{k=-N}^M c_kz^k$ with $c_k\in R(H)$.

Applying Formula~\ref{form:lagrange} to the sum $\sum _{k=-N}^M
c_kz^k$ we get an identity of the form~(\ref{eq:lagrande4}), $f_0$
and $f_1,\dots, f_d$ being polynomials with coefficients in the quotient field of
$R(H)$. Moreover, it is clear from~(\ref{eq:new}) that the
coefficients of $f_0$ are actually in the ring $R(H)$ itself; i.e.
modulo the splitting, $G=T\times H$, $f_0$ is in ~$R(G)$.

To prove Theorem~\ref{thm:lagrange} let
$$
1-a_iz^{k_i}=1-z^{k_i}e^{2\pi\sqrt{-1} \beta_i}=1-e^{2\pi\sqrt{-1}
\alpha_i}
$$
where $\beta_i$ is the restriction of $\alpha_i$ to $H$  and
$k_i=\alpha_i(\xi)$.  Theorem~4 then follows from the formula
above and the identity
\begin{equation}
\label{eq:pi} \pi_i(a_jw^{k_j})=\pi_i
r_i(e^{2\pi\sqrt{-1}\alpha_j}).
\end{equation}

\section{Atiyah-Segal Localization for twisted projective spaces.}
\label{sec:index}

This section describes in detail the twisted projective spaces
which arise as symplectic cuts of stable manifolds. It also
discusses the Atiyah-Segal localization theorem on these twisted
projective spaces.

Let $\alpha_1,\dots,\alpha_m$ be weights of the torus $G$, such
that if $\xi$ is the infinitesimal generator of the circle
subgroup $T$, then $k_i=\alpha_i(\xi)\neq 0$ for $i=1,\dots, m$.
Assume all $k_i$'s are negative.
(At the end of
this section we discuss the case when some $k_i's$ are
positive.)  Let $T$ act on~$\CC^{m+1}$ with weights $k_1, \dots,
k_m, k_{m+1}=-1$. Let $T_\CC$ be the complexification of $T$ which
acts on~$\CC^{m+1}$ with the same weights. The twisted projective
space we are interested in is the orbifold
$$
\widetilde{\CC\PP}^{m}=\CC^{m+1}/\!/
T_\CC=(\CC^{m+1}-\{0\})/T_\CC.
$$
We will not review orbifold theory, but refer the reader
to~\cite{Ruan} for an exposition of orbifold theory and further
references.

To define local orbifold charts on $\tCP^m$, let~$\tilde
{U}_i$ be the $m$-dimensional affine space with coordinates
$(z_1,\dots, \hat z_i,\dots, z_{m+1})$. Denote by $U_i$ the open
subset of $\tCP^m$ where the  projective coordinate~$z_i$ is not
zero. Define the map $\phi_i:\tilde U_i\to U_i$ by
$$
\phi_j(z_1,\dots, \hat z_i,\dots, z_{m+1}) =[z_1,\dots, 1,\dots,
z_{m+1}].
$$
Let $\Gamma_i$ be the finite abelian group of $k_i^{\rm th}$ roots
of unity. Then $w\in \Gamma_i$ acts on $\tilde U_i$ by
$$
w\cdot(z_1,\dots, \hat z_i,\dots, z_{m+1}) = (w^ {k_1}z_1,\dots,
\hat z_i,\dots, w^{k_{m+1}}z_{m+1})
$$
so that $\tilde U_i/\Gamma_i=U_i$. The triples $(\tilde
U_i,\Gamma_i,\phi_i)$ are the orbifold charts of $\tCP^m$.

Assume that $G$ acts on $\CC^{m+1}$ with weights $\alpha_1,\dots,
\alpha_m, \alpha_{m+1}=0$, so that this action descends to an
action on $\tCP^m$. To make sure that this action has a finite
fixed point set we assume that the weights
$\alpha_1,\dots,\alpha_m$ are pairwise linearly independent. Then
the fixed points of the $G$-action on $\tCP^m$ are points
$p_i=[0,\dots,1,\dots,0]$ with $1$ in the $i^\th$ place. An
element $g=\exp(\eta)\in G$ acts on $U_i$ by
\begin{align*}
g\cdot [z_1,\dots, &1,\dots, z_{m+1}]= [e^{2\pi \sqrt{-1}
\alpha_1(\eta)} z_1,\dots, e^{2\pi
\sqrt{-1}\alpha_i(\eta)},\dots,
e^{2\pi \sqrt{-1}\alpha_m(\eta)}z_m,z_{m+1}] \\
&= [e^{2\pi
\sqrt{-1}(\alpha_1(\eta)-\frac{k_1\alpha_i(\eta)}{k_i})}
z_1,\dots,1,\dots, e^{2\pi
\sqrt{-1}(\alpha_m(\eta)-\frac{k_m\alpha_i(\eta)}{k_i})}z_m,
e^{\frac{2\pi \sqrt{-1}\alpha_i(\eta)}{k_i}} z_{m+1}]
\end{align*}
So, the isotropy action of $G$ at $p_i$ is given by the rational
weights $\alpha_j- \frac{k_j\alpha_i}{k_i}$ for $j\neq i$, $1\leq
j\leq m$ and the rational weight $\frac{\alpha_{i}}{k_i}$. Another
way to think about those weights is the following. Notice that if
$G_i$ is the kernel of the map $e^{2\pi\sqrt{-1} \alpha_i}:G\to
S^1$ and $r_i$ is the restriction $R(G)\to R(G_i)$, then
$$
e^{2\pi\sqrt{-1}(\alpha_j-
\frac{k_j\alpha_i}{k_i})}=r_i(e^{2\pi\sqrt{-1}\alpha_j}).
$$
In particular, the rational weights $\alpha_j-
\frac{k_j\alpha_i}{k_i}$ are  genuine integer weights of
$G_i$.  A similar computation shows  that the weights of the action of $G$ at the fixed point $p_{m+1}=[0,\dots,0,1]$ are
$\alpha_1,\dots, \alpha_m$.

We will now explicitly compute the pushforward map in
equivariant K-theory,
$$
\ind_G: K_G(\tCP^m)\to K_G(\pt)
$$
using the Atiyah-Segal localization theorem \cite{AS}; however
before we do this in general we will first consider the situation when all $k_i$'s are equal
to minus one, and $\tCP^m$  is just the standard projective space 
$\CC\PP^m$.

Recall that on a $G$-manifold $X$, the index map (or K-theoretic
pushforward)
$$
\ind_G:K_G(X)\to K_G(\pt)\simeq R(G)
$$
is defined by imbedding $X$ into a linear complex representation
space $V$ of $G$,  applying the Thom isomorphism to map $K_G(X)$ to $K_G(V)$ and then using Bott periodicity to
identify~$K_G(V)$ with $K_G(\pt)$.  We also recall that by the Atiyah-Segal localization theorem the
restriction map
$$
K_G(X)\to K_G(X^G)
$$
becomes an isomorphism after localizing with respect to a certain
prime ideal of $R(G)$ (see~\cite{AS}), and that  it is possible to
write an explicit formula for $\ind_G(\delta)$ in terms of this
restriction. In the case $X=\CC\PP^m$ this formula says that
\begin{equation}
\label{eq:torsion-free} \ind_G (\delta)=\sum_{i=1}^m
\frac{\delta_i }{(1-e^{2\pi \sqrt{-1}\alpha_i})\prod_{j\neq i}
(1-e^{2\pi \sqrt{-1} (\alpha_{j}-\alpha_i)})} +
\frac{\delta_{m+1}}{\prod_{j=1}^m (1-e^{2\pi
\sqrt{-1}\alpha_{j}})},
\end{equation}
where $\delta_i$ is the restriction of $\delta$ to $p_i$, and the
denominators in the formula are just the virtual characters of the
exterior algebra complexes $\sum (-1)^k\Lambda^kT_{p_i}$ of the tangent
spaces
$T_{p_i}$ at the fixed points.

We now drop the assumption that all $k_j$'s are equal to minus one, so
that $\tCP^m$ may have orbifold singularities. Let us recall the
definition of the equivariant index map for orbifolds (for more
details see~\cite{Kawasaki} or~\cite{Paradan}). Assume an orbifold
$X$ is presented as a quotient $Y/K$ of a manifold~$Y$ by a locally
free action of a compact group $K$ action. (We will describe such a
presentation of $\tCP^m$ shortly.) Moreover, assume $G$ acts on
$Y$ and commutes with the $K$ action, then the action of $G$
descends to an action on $X$. It is well-known there exists an
isomorphism, (which for purposes of this paper will be treated as a
definition of $K_G(X)$)
$$
\Psi:K_G(X) \stackrel{\simeq}{\longrightarrow} K_{G\times K} (Y).
$$
Then for $\delta \in K_G(X)$ we define its index by
$$
\ind^X_G\delta =\Big(\ind^Y_{G\times K}
\big(\Psi(\delta)\big)\Big)^K\in R(G),
$$
the $K$ invariant part of the $G\times K$ index of
$\Psi(\delta)$.  The relative version of the localization theorem  on~$Y$
produces a localization formula for orbifolds, which expresses
$\ind^X_G(\delta)$ in terms of the restriction of~$\delta$ to
$X^G$. We will not give the general version of this formula, since
we will only apply it to the case of the weighted projective space
$\tCP^m$. So from this point on we specialize to the case $X=\tCP^m$.

The twisted projective space $\tCP^m$ can be realized as a the
symplectic reduction of $\CC^{m+1}$ by the action of $T$, which
is just the quotient $S^{2m+1}/T$, where $S^{2m+1}$ is the sphere
$$
S^{2m+1} =\{(z_1,\dots, z_{m+1})\in \CC^{m+1} |\sum
k_i|z_i|^2=-1\}
$$
on which the circle $T$ acts locally freely. So, as mentioned
above there is an isomorphism,
$$
\Psi:K_{G}(\tCP^m)\stackrel{\simeq}{\longrightarrow}K_{G\times
T}(S^{2m+1}).
$$

For $\delta\in K_G(\tCP^m)$ we are interested in computing the
 K-theoretic index
$$
\ind_G^{\tCP^m} \delta= \Big(\ind^{S^{2m+1}}_{G\times T}
\Psi(\delta)\Big)^{T}\in R(G)
$$
by means of  Atiyah-Segal localization.  Recall that the fixed
points of the $G$ action on $\tCP^m$ are the
points $p_1,\dots, p_{m+1}$.  Denote by $\iota_i$ the inclusion
$p_i\to \tCP^m$. Let  $s_i$ be the circle inside $S^{2m+1}$
which after dividing by $T$ becomes $p_i$ and let
$\tilde\iota_{i}:s_i\to S^{2m+1}$ be the natural inclusion.  The stabilizer of the $G\times T$ action on $s_i$ is the group
$G_i\times T'$, where $T'$ is the subgroup \mbox{$\{(t,t^{-1})|t\in
T\}$} of $G\times T$ and~$G_i$ is the kernel of
$e^{2\pi\sqrt{-1}\alpha_i}: G\to S^1$. Then
$$
K_{G\times T}(s_i)=K_{G_i\times T'} (p_i)\cong R(G_i\times T').
$$

We now define  maps $\tilde r_i$ and $\tilde \pi_i$, which,
as we will explain below, are extensions of the maps,~$r_i$
and~$\pi_i$ of Section~\ref{sec:lagrance}. The map
$\tilde r_i$ will be the restriction homomorphism
\begin{equation}
\label{eq:r} R(G\times T')\to R(G_i\times T'),
\end{equation}
which clearly extends $r_i$, the restriction map from $R(G)$ to
$R(G_i)$.

Consider the covering map $\rho:G_i\times T'\to G$ which sends
$(g,(t,t^{-1}))\in G\times T'$ to $gt\in G$. The kernel of this
map is the group $\Gamma_i\cong G_i\cap T$ which can be identified with
the $k_i^{\rm th}$ roots of unity, so that
$$
G_i\times T'/\Gamma_i= G.
$$

The Gysin map
\begin{equation}
\label{eq:pi2} \tilde \pi_i: K_{G\times T}(s_i)=K_{G_i\times
T'}(\pt)=R(G_i\times T') \to R(G)=K_G(\pt)=K_G(p_i)
\end{equation}
is defined as follows. Let $V$ be a virtual representation for
$G_i\times T'$, and let $\tilde \pi_i(V)$  be the subspace
$V^{\Gamma_i}$ of $V$ fixed by $\Gamma_i$. Since $G=G_i\times
T'/\Gamma_i$, the space $V^{\Gamma_i}$ is a virtual representation
space for $G$, so $\tilde \pi_i$ is well-defined. Again this map is 
clearly an extension of the map $\pi_i:K_{G_i}(\pt) \to K_H(\pt)$
defined in~(\ref{eq:pi}), where $H$ is a subtorus for which
$G=H\times T$.

The key ingredient in the Atiyah-Segal localization formula is the
fact that the composition of the pushforward map
$\tilde\iota_{i!}$ and the restriction $\tilde \iota^*$
$$
K_{G_i\times T'}(\pt)=K_{G\times
T}(s_i)\stackrel{\tilde\iota_{i!}}{\to} K_{G\times T} (S^{2m+1})
\stackrel{\tilde \iota^*_{i}}{\to}K_{G\times T}(s_i)=K_{G_i\times
T'}(\pt)
$$
is multiplication by the character of the exterior algebra complex
$\sum(-1)^k\Lambda^kT_{p_i}$ of the tangent space at $p_i$, which is just
\begin{equation}
\label{eq:multip} \tilde r_i\Big((1-\zeta)\prod_{j\neq i}
\big(1-e^{2\pi \sqrt{-1} \alpha_{j}}\big)\Big),
\end{equation}
where $\zeta$ is the generator of the character
ring $R(T')$. (Recall that $G$ acts on $T_{p_i}$ with the rational
weights, but the group $G_i\times T'$, a cover of $G$, acts on
$T_{p_i}$ with integer weights. So, in~(\ref{eq:multip}) $\zeta$
is the character of the representation of $T'$ associated with the
rational weight $\frac{\alpha_i}{k_i}$ of $G$.)

The Atiyah-Segal localization theorem states that the map
$\tilde\iota_!=\sum \tilde\iota_{i!}$ becomes an isomorphism after
 localization. Thus, together with the fact that $\tilde
\iota^*\tilde\iota_{i!}$ is multiplication by~(\ref{eq:multip}) it
implies that
$$
\ind_{G\times T}^{S^{2m+1}} (\Psi(\delta))=\sum_{i=1}^m
\frac{\tilde\iota^*_i \Psi(\delta)}{\tilde
r_i\Big((1-\zeta)\prod_{j\neq i} \big(1-e^{2\pi \sqrt{-1}
\alpha_{j}}\big)\Big)} +
\frac{\tilde\iota^*_{m+1}\Psi(\delta)}{\prod_{i=1}^m (1-e^{2\pi
\sqrt{-1}\alpha_{i}})}.
$$
It remains to take $T$ invariants of both part of this formula.
Notice that taking the $T$-invariant part of a K-class in 
$K_{G\times T}(s_i)=K_{G_i\times T'}(\pt)$ is equivalent to
applying the map $\tilde \pi_i$ to this class. Hence we get

\begin{formula}
\label{form:index} For $\delta\in K_G(\tCP^m)$
\begin{equation}
\label{eq:ASCP} \ind_G^{\tCP^m} (\delta)=\sum_{i=1}^m \tilde
\pi_i\Big(\frac{\tilde \iota^*_i \Psi(\delta)}{\tilde
r_i\Big((1-\zeta)\prod_{j\neq i} \big(1-e^{2\pi \sqrt{-1}
\alpha_{j}}\big)\Big)}\Big) + \frac{\tilde
\iota^*_{m+1}\Psi(\delta)}{\prod_{i=1}^m (1-e^{2\pi
\sqrt{-1}\alpha_{i}})}.
\end{equation}
\end{formula}

In  case  not all the numbers $k_i$ are
negative, assume that for the first $r$ weights,
$\alpha_1,\dots,\alpha_r$, these numbers  are positive and for the
others are negative. Then apply
Formula~\ref{form:index} to the twisted projective space defined
for the weights
$-\alpha_1,\dots,-\alpha_r,\alpha_{r+1},\dots,\alpha_m$ and the
equivariant K-theory class~$\delta\prod_{i=1}^r
e^{-2\pi\sqrt{-1}\alpha_i}$. This, after an easy computation,
yields

\begin{formula}
\label{form:index-sign} For $\delta\in K_G(\tCP^m)$
\begin{align*}
(-1)^r\ind_G^{\tCP^m} (\delta\prod_{i=1}^r
e^{-2\pi\sqrt{-1}\alpha_i})=&\sum_{i=1}^m \sgn(-k_i)\tilde
\pi_i\Big(\frac{\tilde \iota^*_i \Psi(\delta )}{\tilde
r_i\Big((1-\zeta)\prod_{j\neq
i} \big(1-e^{2\pi \sqrt{-1} \alpha_{j}}\big)\Big)}\Big) \\
 &+\frac{\tilde \iota^*_{m+1}\Psi(\delta)}{\prod_{i=1}^m
(1-e^{2\pi \sqrt{-1}\alpha_{i}})}.
\end{align*}
\end{formula}

\section{Calculation of local indices for GKM spaces.}
\label{sec:local-index}

In this section we apply Formula~\ref{form:index} to the
twisted projective spaces, $S_\e$, to prove Lemma~\ref{lem:local} and
Theorem~\ref{thm:GKM}. We also describe the relationship between
formulas~(\ref{eq:lagrande-intro}) and~(\ref{eq:ASCP}).

Let us recall the definition of symplectic cuts of stable
manifolds. Assume $(M,\omega)$ is a Hamiltonian $G$-space. Choose
a generic circle subgroup $T$ of $G$ such that $M^T=M^G$. Let $p$
be an isolated fixed point  and $S$ the stable
manifold of $p$. Let $\omega_S$ be the restriction of $\omega$ to
$S$, and consider the space $S\times \CC$ with the symplectic form
$(w_S,-\sqrt{-1}dy\wedge d\bar y)$, where $y$ is the
complex coordinate on $\CC$. If $T$ acts on $\CC$ with the weight
$-1$, the action of $T$ on $S\times \CC$ is Hamiltonian. Restrict
the moment map~$\phi$ of the $T$ action to $S$. If $\phi(p)=c$,
 the moment map for the action of $T$ on $S \times \CC$ is
$$
\tilde \phi(x,y)=\phi(x)-c-|y|^2 \, ,
$$
and the symplectic cut $S_\e$ is the symplectic
reduction of $S\times \CC$ at $-\e$. For more details and an
explanation of why this definition of $S_\e$ coincides with
definition~(\ref{eq:cut}) see~\cite{Lerman}.

The Kirwan  map
$$
\tilde \kappa_\e:K_{G\times T}(S\times \CC)\to K_G(S_\e)
$$
is the composition of the restriction of a class in $K_{G\times
T}(S\times\CC)$ to  $\tilde \phi^{-1}(\e)$ and the identification
of $K_{G\times T}(\tilde \phi^{-1}(\e))$ with $K_G(S_\e)$.  As before, let $T'$ be the circle subgroup of $G\times T$ given by
$\{(t,t^{-1})|t\in T\}$, so that there is a canonical
identification $G\times T\cong G\times T'$. Notice that $T'$ act
trivially on $S$, hence
$$
K_{G\times T}(S\times \CC)=K_{G\times T'}(S\times
\CC)=K_{G}(S)\otimes K_{T'}(\CC).
$$
Using this identification  define, for a class, $a \in K_G (M)$
$$
\kappa_\e (a)=\tilde \kappa_\e(a_S\otimes 1),
$$
where $a_S$ is the restriction of $a$ to $S$.
Then $I_p(a)$ is just the $G$-equivariant index of $\tilde
\kappa_\e(a_S\otimes 1)$.

\begin{remark} \label{rem:explanation} We will explain why $I_p$ is an 
  $R(G/T)$ (rather than $R(G)$) module homomorphism. There are
  two $R(G)$ module structures on $K_{G\times T}(S\times\CC)$:
  one coming from the projection of $G\times T$ onto the first
  factor and the other coming from the multiplication map
  $G\times T\to G$, which maps $(g,t)\mapsto g\cdot t^{-1}$, and
  the map $\tilde \kappa_\e$ is an $R(G)$-module homomorphism
  with respect to the second (not the first) $R(G)$-module
  structure.  However for $I_p$ to be an $R(G)$-module
  homomorphism $\tilde \kappa_\e$ has to be an $R(G)$-module
  homomorphism with respect to the first module structure.
  Because of this~$I_p$ is just an $R(G/T)$-module homomorphism.
  (This is also clear, by the way, from the algebraic description
  of the map, $I_p$, in Theorem~4.  Namely formulas~(3.4) and
  (3.6) make clear that the map, $f \to f_0$ in Theorem~4 is
  \emph{not} an $R(G)$ morphism.)
\end{remark}

Let~$T_p$ be the tangent space at $p$. This space decomposes into
a direct
sum
$$
T_p=T_p^-\oplus T_p^+,
$$
where $T$ acts on $T^\pm_p$ with positive and negative weights
respectively. The exponential map identifies a neighborhood of $p$
in $S$ with a neighborhood of the origin inside $T_p^-$. Let $G$
act with weights $\alpha_1,\dots, \alpha_m$ (so far we allow
linear dependencies) on $T_p^-$. Then, it is obvious that $S_\e$
can be identified with the projective space $\tCP^m$ defined in
the previous section.

Before we specialize to GKM spaces we will prove
Lemma~\ref{lem:local}.

\begin{proofof}{Lemma~\ref{lem:local}}
As above identify the symplectic cut of the stable
manifold $S$ at~$p$ with $\tCP^m$. Assume for the moment that we are 
in the GKM setting so that every
pair of $\alpha_i$'s are linearly independent. Apply
Formula~\ref{form:index} to $\kappa_\e(\tau)$, where $\tau$ is a
class with $\tau(q)=0$ whenever $q<p$. In particular if a descending 
edge goes from $p$ to $q$, then $\tau(q)=0$. Hence $\tilde \iota^*_i\Psi(
\kappa_\e(\tau))=0$, unless $i=m+1$ and
$$
\tilde \iota^*_{m+1}\Psi(\kappa_\e(\tau))=\tau_p.
$$
Thus Formula~\ref{form:index} yields
$$
I_p(\tau)=\ind_G^{\tCP^m} (\kappa_\e(\tau))= \frac{\tilde
\iota^*_{m+1}\Psi(\kappa_\e(\tau))}{\prod_{i=1}^m (1-e^{2\pi
\sqrt{-1}\alpha_{i}})}=\frac{\tau_p}{\Lambda^-_p}
$$
which finishes the proof in  the case $\alpha_i$'s are pairwise
linearly independent.

If we allow linear dependencies among pairs of $\alpha_i$'s,
then a very similar proof will go through. Atiyah-Segal
localization gives a formula analogous to
Formula~\ref{form:index}, the difference being that
the fixed point set of the $G$ action on $\tCP^m$ contains an
isolated point $p_{m+1}$ and other (possibly nonisolated) fixed
points.   Here the first term on the right hand side of (\ref{eq:ASCP}) 
(which is the
contribution to the index of these other points) is more
complicated.  However if the class $\tau$ is supported above $p$, that is 
$\tau(q)=0$ whenever $\phi(q)<\phi(p)$, then it is clear that the 
first term on the right hand side of~(\ref{eq:ASCP}) vanishes, since 
restrictions to 
the corresponding fixed point sets are zeros.

To finish the proof we will explain why we can
assume  without loss of generality that $\phi(q)<\phi(p)$ if and only if $q<p$, so that the condition
that $\tau(q)=0$ for  $q<p$ implies that $\tau$ is supported above
$p$. This follows from the equivariant version of the
following standard result in classical Morse theory (see
\cite[Theorem~4.1]{Milnor}). Given two real 
numbers
$a<b$ assume $\phi^{-1}([a,b])$ contains two critical points~$p$ and~$q$
whose stable and unstable manifolds do not intersect. Then for any
$a<c,c'<b$, there exists another Morse function $\phi'$ which coincides
with~$\phi$ outside of $\phi^{-1}([a,b])$, has the same critical points 
as
$\phi$ and the same stable and unstable manifolds as $\phi$, and
$\phi'(p)=c$ and $\phi'(q)=c'$. It is clear from the proof of
\cite[Theorem~4.1]{Milnor} that the proof of this fact does not rely on
any transversality arguments and can be restated and proved in the
equivariant setting.  \end{proofof}

Let us now turn to the proof of Theorem~\ref{thm:GKM}. We assume
that $M$ is a GKM manifold with respect to the $G$ action, so
that, the one-skeleton, $\Gamma$ with its $G$-action carries all the information we need to compute
the equivariant cohomology and the equivariant K-theory of~$M$.

Let $p\in M$ be a fixed point under the $G$ action, or, in other
words, a vertex of $\Gamma$, and let $T$ act on $T_p^-$ with the
weights $\alpha_1,\dots, \alpha_m$. The GKM assumption implies
that every pair of weights is linearly independent, and in  terms of the
graph $\Gamma$ it means that the descending edges $e_1,\dots, e_m$
coming out of $p$ are labeled by the weights $\alpha_1,\dots,
\alpha_m$.

For any class $a\in K_G(M)$, let us apply
Formula~\ref{form:index} to the class $\kappa_{\epsilon}(a)$ 
to compute $I_p(a)$.  We can make Formula~\ref{form:index} more specific
by computing the restrictions of $\kappa_\e(a)$ to the fixed
points of $\tCP^m$. It is easy to see that
$$
\tilde \iota_i^* \Psi(\kappa_\e(a))=r_i(a_p)=r_i(a_{q_i}),
$$
where $q_i$ is the second vertex of the edge $e_i$.  By inserting
these identities in (\ref{eq:ASCP}) one immediately gets the 
identity~(\ref{eq:GKM})
of Theorem~\ref{thm:GKM}.

We also observe that if we  apply the Lagrange interpolation
formula~(\ref{eq:lagrande-intro}) to the virtual character $a_p\in
R(G)$ with $\alpha_1,\dots,\alpha_m$ being the labels of the edges
of $\Gamma$ pointing down from $p$, we get a slightly
weaker version of Theorem~\ref{thm:GKM}. Namely, we get a
formula whose right hand side is identical to the right hand
side of~(\ref{eq:GKM}), but has the term $f_0$ (instead
of $I_p(a)$) on the left hand side. So, the proof of
Theorem~\ref{thm:GKM} presented above gives a geometric
interpretation of the term $f_0$ in the Lagrange interpolation
formula. 
Similarly,
Formula~\ref{form:index-sign} can be interpreted as the geometric
analogue of the formula~(\ref{eq:lagrande-intro}).

Notice also that the above discussion can be adapted to the
setting of
equivariant cohomology, with the pushforward in K-theory, or
index, replaced by the pushforward in cohomology, or
integration.

\section{The restriction of the classes $\tau_p$ to fixed points.}
\label{sec:restriction} Theorem~\ref{thm:GKM} enables one to express local indices in terms of restrictions of
equivariant K-classes to fixed points. In this
section we will do the opposite. Namely, we will 
compute the restrictions of equivariant K-theory classes to fixed
points in terms of local indices. Since the $\tau_p$'s generate
$K_G(M)$ it is enough to do this for the classes
$\tau_p$'s. Hence, we will derive an explicit graph-theoretic formula
for the restriction of  $\tau_p$ to a fixed point $q\in M^G$.

To state this result we recall some notation: Given an oriented
edge $e$ of $\Gamma$ we denote by $i(e)$ and~$t(e)$ the initial
and terminal vertices of $e$, and we denote by $\bar e$ the edge
obtained from $e$ by reversing its orientation. Thus $i(\bar
e)=t(e)$ and $t(\bar e)=i(e)$. We call $e$ an ascending edge if
$\alpha_e(\xi)>0$ and a descending edge if $\alpha_e(\xi)<0$.

Let $H$ be a complimentary torus to $T$ in $G$. As in
Section~\ref{sec:lagrance} we will make use of the splitting
$G=T\times H$ to identify the ring $R(G)$ with the ring of finite
sums
\begin{equation}
\label{eq:R} \sum_{k=-N}^M c_kz^k, \ \  c_k\in R(H),
\end{equation}
and we will denote by $Q(H)$ the quotient field of the ring $R(H)$
and by $\hat R(G)$ the ring of finite sums
\begin{equation}
\label{eq:Q} \sum_{k=-N}^M c_kz^k, \ \  c_k\in Q(H).
\end{equation}
Thus by~(\ref{eq:R}) $\hat R(G)$ contains $R(G)$ as a subring.

Let $e$ be an ascending edge of $\Gamma$. The key ingredient in
our combinatorial formula for $\tau_p(q)$ is a $Q(H)$ module
endomorphism
$$
Q_e:\hat R(G)\to  \hat R(G).
$$

To define this endomorphism, let $p=t(e)$ and let $e_1,\dots
,e_{r+1}$ be the descending edges of $\Gamma$ with $i(e_j)=p$. We
will order these edges so that $e_{r+1}=\bar e$. Let
$\alpha_j=\alpha_{e_j}$. Then
$e^{2\pi\sqrt{-1}\alpha_j}=z^{k_j}e^{2\pi\sqrt{-1}\beta_j}$, where
$k_j$ is a negative integer and $\beta_j$ is an element of the
weight lattice of $H$. Denote by $G_e$ the kernel of
$e^{2\pi\sqrt{-1}\alpha_e}:G\to S^1$. Given an element $f$ of
$\hat R(H)$ and $h\in H$ we define $Q_e f$ to be the expression
\begin{align}
\label{eq:Q_e} Q_e f(z,h)&=\Big(
\prod_{j=1}^{r+1}(1-z^{k_j}e^{2\pi\sqrt{-1}\beta_j}) \Big) \pi_i
r_i\,\, \frac{f(w,h)}{(1-\frac{z}{w})\prod_{j=1}^r(1-
e^{2\pi\sqrt{-1}\alpha_j})}\\
\noalign{\hbox{or alternatively the sum}}
 Q_e f(z,h) &=\Big( \prod_{j=1}^{r+1}(1-z^{k_j}e^{2\pi\sqrt{-1}\beta_j}) \Big)
\frac{1}{k} \sum_{i=1}^r
\frac{f(w_i,h)}{(1-\frac{z}{w_i})\prod_{j=1}^r(1-w_i^{k_j}
e^{2\pi\sqrt{-1}\beta_j})},
\end{align}
where $w_1,\dots, w_k$ are the preimages of $h$ in $G_e$ with
respect to the projection $G_e \to H$, and $r_i$ and $\pi_i$ are
the restriction and Gysin maps. By Theorem~\ref{thm:lagrange},
this expression is a finite sum of the form~(\ref{eq:Q}) and hence
an element of $\hat R(G)$.

Notice that with $e=e_i$, $Q_e f(z,h)$ is the term
$f_i$ in formula~(\ref{eq:lagrange-intro2}) of 
Theorem~\ref{thm:lagrange} multiplied by $
\prod_{j=1}^{r+1}(1-z^{k_j}e^{2\pi\sqrt{-1}\beta_j}) $. From the
geometric interpretation of Theorem~\ref{thm:lagrange}, as 
Atiyah-Segal localization on twisted projective spaces, it is
possible to interpret $Q_e$ as a purely topological operation. This is an easy
exercise, and we omit the details.

Now let $\gamma$ be a path in $\Gamma$ joining $p$ to $q$; i.e. a
sequence of oriented edges $e_1,\dots, e_s$ with $i(e_1)=p_1$,
$t(e_s)=q$ and $t(e_j)=i(e_{j+1})$. We will call $\gamma$ an
ascending path if all the $e_j$'s are ascending, and we will
denote by
$$
Q_\gamma: \hat R(G)\to \hat R(G)
$$
the composition $Q_{e_s}Q_{e_{s-1}}\dots Q_{e_1}$.

\begin{theorem}
\label{thm:paths} The restriction $\tau_p(q)$ of $\tau_p$ to $q\in
M^G$ is equal to the sum
\begin{equation}
\label{eq:sum} \sum Q_\gamma(1)
\end{equation}
over all ascending paths $\gamma$ in $\Gamma$ joining $p$ to $q$.
\end{theorem}

\begin{proof} For $p=q$ this formula is obvious, so we assume
$q\neq p$. Let $e_1,\dots , e_r$ be the ascending edges of
$\Gamma$ with terminal vertex $q$ and let $q_j=i(e_j)$. Since
$I_q(\tau_p)=0$ we get from Theorems~\ref{thm:GKM}
and~\ref{thm:lagrange}
\begin{equation}
\label{eq:sum2} \tau_p(q)=\sum Q_{e_i}\tau_p(q_i).
\end{equation}
If one of the $q_i$'s, say $q_1$ is equal to $p$, one can incorporate it
as a term in the sum~(\ref{eq:sum}), since $e_1$ is a path of
length one joining $p$ to $q$. For the $q_i$'s which are not equal
to $p$, let $e_{ij}$ be the ascending edges of $\Gamma$ with
$t(e_{ij})=q_i$ and let $q_{ij}=i(e_{ij})$. By
iterating~(\ref{eq:sum2}) we get for the contribution of these~$q_i$'s to 
the sum~(\ref{eq:sum2})
\begin{equation}
\label{eq:sum3} \sum_{i,j}Q_{e_i}Q_{e_{ij}}\tau_p(q_{ij}).
\end{equation}
Again, if one of the $q_{ij}$'s is equal to $p$, the summand
in~(\ref{eq:sum3}) corresponding to it gets incorporated in the
sum~(\ref{eq:sum}), since the path with edges $e_{ij}$ and $e_i$
is an ascending path on $\Gamma$ joining $p$ to $q$. For the
remaining $q_{i,j}$'s we iterate this argument again. It is clear
that after sufficiently many iterations we obtain the
formula~(\ref{eq:sum}).
\end{proof}

An analogue of Theorem~\ref{thm:paths} for the equivariant cohomology ring
$H^*_G(M)$ can be found in~\cite{GZ3}. In equivariant cohomology, the
$\tau_p$'s are defined as the equivariant Thom classes of the unstable
manifolds, and it was shown in~\cite{GZ3} that the restriction of these
classes to the fixed points are given by a formula analogous
to~(\ref{eq:sum}) with operators, $Q_e$, defined by formulas similar
to~(\ref{eq:Q_e}).  Moreover, Catalin Zara was able to show that, for this
combinatorial version of~(\ref{eq:sum}) a lot of summands cancel each
other out making this formula an effective tool for computational
purposes.  We conjecture (or at least hope) that the same will be true in
K-theory.

\section{Example: the case of Grassmannian.}
\label{sec:example}

In this section we present a generalization of
the definition of local index. This allows us to single out a
class of GKM manifolds for which the computation of the
restriction of the class, $\tau_p$, to $M^G$ is not much more
complicated in $K$-cohomology than in ordinary cohomology.  (We
will show that one example of such a manifold is the Grassmannian.)

For a compact symplectic manifold $(M, \omega)$ with a Hamiltonian
$G$ action the total index map $\mathcal I : K_G(M)\to K_G(M^G)$
depends on a choice of a circle subgroup $T\in G$ with $M^G=M^T$.
Recall that the local index of $a\in K_G(M)$ at a fixed point $p$
depends on the symplectic cut of the stable manifold $S_p$ of $\phi$ at
$p$
with respect to the $T$ action. We now modify this definition of
local index by considering  symplectic cuts of these stable manifolds with 
respect to other circle subgroups of $G$.

Namely, besides choosing a circle subgroup $T$ of $G$ with
$M^G=M^T$, choose for each fixed point~$p$ a circle subgroup $T_p$
of $G$, such that the $T_p$-moment map $\mu_p: M\to \RR$
restricted to $S_p$ attains its maximum at $p$. This allows us to
define this local symplectic cut using the $T_p$ action instead of
the~$T$ action. So, now let $S_\e$ be the symplectic cut of
the stable manifold of $\phi$ at $p$ with respect to the circle~$T_p$
action. As before we have a map $\kappa_\e:K_G(M)\to K_G(S_\e)$
and we define the local index $\tilde I_p$ by
$$
\tilde I_p(a)=\ind_G(\kappa_\e(a)).
$$

Notice that all the results in Section~\ref{sec:proofs} hold  for
the new local indices $\tilde I_p$ . In particular,
Theorem~\ref{thm:main} is still true. (Notice, however, that the
classes $\tau_p$ may be different from those defined using the
old definition of local indices.) Also Formulas~\ref{form:index}, 
\ref{form:index-sign} as well as Theorem~\ref{thm:GKM} hold with the 
following minor changes.
If $\xi_p$ is the infinitesimal generator of $T_p$, then the
numbers~$k_i$ have to be defined as~$\alpha_i(\xi_p)$, the
maps $\tilde \pi_i $ and $\tilde r_i$ are defined using the circle
$T_p$ instead of $T$, and $\zeta$ is the generator of the character ring 
$R(T_p)$.

We now define the new total index $\tilde {\mathcal I}:K_G(M)\to
K_G(M^G)$ to be the sum of all $\tilde I_p$. (Notice that the new total 
index is no longer an $R(G/T)$-module homomorphism.) Let us call
$\tilde {\mathcal I}$  \emph{torsion-free} if we can pick the circles
$T$ and $T_p$'s in such a way that all the numbers $k_i$'s for all
the fixed points are equal to minus one. Then Theorem~{\ref{thm:GKM}
simplifies to

\begin{theorem}
\label{thm:torsion-free} Let $M$ be a GKM space.  For $p\in
V=M^G$, let  $e_1,\dots, e_m$ be the descending edges
with initial vertex  at $p$. Let~$e_i$ connect $p$ to $q_i$ and be
labeled by weight $\alpha_i$. If $\tilde {\mathcal I}$ is
torsion-free, then for any $a\in K_G(M)$ we have
\begin{equation}
\label{eq:torsion-free2} \tilde I_p(a)=\sum_{i=1}^m
\frac{r_i(a_{q_i})}{(1-e^{2\pi\sqrt{-1}\alpha_i})\prod_{j\neq i}
\big(1-e^{2\pi \sqrt{-1} (\alpha_{j}-\alpha_i)}\big)} +
\frac{a_p}{\prod_{i=1}^m (1-e^{2\pi \sqrt{-1}\alpha_{i}})},
\end{equation}
\end{theorem}

\begin{proof} This  immediately follows from~(\ref{eq:GKM}) with all $k_i$'s equal to
minus one, or from a straightforward application
of~(\ref{eq:torsion-free}).
\end{proof}

As mentioned above local indices can be defined in the setting of 
equivariant cohomology by using integration instead of the index map for 
the pushforward. Then Theorem~\ref{thm:torsion-free} has a counterpart in 
the equivariant cohomology which states that for $a\in H^*_G(M)$,
\begin{equation}
\label{eq:torsion-free3} \tilde I_p(a)=\sum_{i=1}^m
\frac{r_i(a_{q_i})}{\alpha_i\prod_{j\neq i}
\big(\alpha_{j}-\alpha_i\big)} +
\frac{a_p}{\prod_{i=1}^m \alpha_{i}},
\end{equation}
where $\tilde I_p$ is the local index in equivariant cohomology.
Notice that~(\ref{eq:torsion-free3}) is obtained from~(\ref{eq:torsion-free2}) by the formal substitution of every expression 
of the form $(1-e^{2\pi\sqrt{-1}\alpha})$ by~$\alpha$. 

Without giving  details we remark that it is possible to
prove an analogue of Theorem~\ref{thm:paths} for this new local index. However the operators $Q_\gamma$ will be defined
differently and will no longer map $\hat R(G)$ into itself, where
$\hat R(G)$ is the ring of the sums of the form~(\ref{eq:Q}). In
the  torsion free case
 the operators $Q_\gamma$ will simplify 
and after
formal substitution of $(1-e^{2\pi\sqrt{-1}\alpha})$ by~$\alpha$ will
resemble the operators of~\cite{GZ3} appearing in their ``path integral
formula''.

We will conclude this section by showing that the Grassmannian,
$Gr (k,n)$, of $k$-planes in $\CC^n$ is an example of a space for
which one can define a torsion free total index map.  To see
this, identify  the $n$-dimensional torus  $\tilde G$ with the product of
$n$ circles $S_1\times \dots \times S_n$. Let $S_i$ act on the
$i^{\th}$ component of $\CC^n$ with weight  $1$ and with weight
$0$ on the other components. This action induces an action of
$\tilde G$ on $Gr(k,n)$. If  $S$ is the diagonal of the
torus $\tilde G$ then its action on $Gr(k,n)$ is trivial, and the
action of $G=\tilde G/S$ on $Gr(k,n)$ is effective. The $G$ action
on $Gr(k,n)$ is known to be GKM and its one-skeleton 
$\Gamma$ is the Johnson graph.

Let $\xi_1,\dots, \xi_n$ be the infinitesimal generators of the
circles $S_1,\dots, S_n$ . They form a basis of the Lie algebra
$\tilde { \mathfrak g}$ of $\tilde G$. Let
$\alpha_1,\dots,\alpha_n$ be the dual basis of $\tilde { \mathfrak
g}^*$. Then $\sum c_i\alpha_i$ with $c_i\in \ZZ$ is a weight
of~$G$ as long as $\sum c_i=0$.

Let $v_i$ be a nonzero vector in the $i^{\rm th}$ component of
$\CC^n=\CC\times\dots\times \CC$. The fixed points of the~$G$
action on $Gr(k,n)$ are indexed by the $k$-element subset of
$\{1,\dots,n\}$. Namely if $I$ is such a set, the fixed point
$p_I$ is the span of vectors $v_{i_k}$ with $i_k\in I$. The
weights of the isotropy action at $p_I$ are $\alpha_i-\alpha_j$
with $i\notin I$ and $j\in I$. Let $\tilde T_I$ be the diagonal 
subcircle of the torus $\prod_{i\in I} S_{i}$ and $T_I$ its image 
inside $(\prod_{i\in I} S_{i})/S$. Then the isotropy action of
$T_I$ at $p_I$ is given by the weight $-1$, and the moment map associated 
to the $T_I$ action attains its maximum at $p_I$.

So, if we pick any generic circle subgroup $T$ of $G$ and then
define local indices $\tilde I_{p_I}$ using symplectic cuts with
respect to the actions of $T_I$, the total index will be
torsion-free. In particular, Theorem~\ref{thm:torsion-free} applies to 
$Gr(k,n)$.

Now let us specialize to the case of the ordinary (nonequivariant)
cohomology and K-cohomology rings of $Gr(k,n)$. They are known to be
isomorphic, where the isomorphism $\Phi:K(Gr(k,n))\to
H^*(Gr(k,n))$ is given by sending the Chern classes in K-theory of
the dual of the tautological vector bundle on $Gr(k,n)$ to the
corresponding Chern classes in cohomology. For more details concerning
this isomorphism and the discussion below see~\cite{Lenart}.

Pick a generic circle, $T$, in $G$, such that its moment map  attains its minimum at $p_I$ with 
$I=\{1,\dots,k\}$. Call the
closure of the stable manifold at $p_I$ the Schubert variety~$X_I$.  Then
we can define two different bases of $H^*(Gr(k,n))$  one basis,~$s_I$, being given by the Poincare duals, or
Thom classes, of the $X_I$'s, and  the other basis, $g_I$, being given by the images
under $\Phi$ of the structure sheaves of the $X_I$'s. (The topological
K-theory of $Gr(k,n)$ can be identified with its algebraic K-theory, so
that we can use coherent sheaves, in particular the structure sheaves of
Schubert varieties, to define K-classes, see~\cite{BFM} for details.) 
These
bases are not the same and the transition matrix between these two bases
was worked out in \cite{Lenart} using the combinatorics of Schur and
Grothendieck polynomials.

Notice that the classes $\tau_{p_I}$ constructed in Theorem~\ref{thm:main} 
descend
to a basis $\hat\tau_{p_I}$ of $K(Gr(k,n))$. From the discussion
above it is clear that $\Phi(\hat \tau_{p_I})=s_I$. So this allows one to
interpret the coefficients computed in~\cite{Lenart}
geometrically. Namely these are the coefficients appearing when we express 
K-theory classes of the structure sheaves of the $X_I$'s as linear 
combinations 
of the $\hat \tau_p$'s. Hence, the problem of computing
these coefficients  can be reduced to the computation of the
(nonequivariant) local indices of the K-theory classes of the
structure sheaves of the $X_I$'s. It would be very interesting to
reprove the results of~\cite{Lenart} using this geometric
approach.

\end{document}